\documentclass[11pt,twoside]{article}
\usepackage{rotating}
\usepackage{amsmath}
\usepackage{amsfonts}
\usepackage{amssymb}

\newtheorem{thm}{Theorem}[section]
 
 \newtheorem{lem}[thm]{Lemma}
  \newtheorem{defi}[thm]{Definition}

\newcommand{\brdef}{\begin{defi}}
\newcommand{\erdef}{\end{defi}}

\newcommand{\blem}{\begin{lem}}
\newcommand{\elem}{\end{lem}}

\newtheorem{corr}{Corollary}[section]
\newcommand{\bcor}{\begin{corr}}
\newcommand{\ecor}{\end{corr}}

\newtheorem{key}{keyword}[section]
\newcommand{\bkey}{\begin{key}}
\newcommand{\ekey}{\end{key}}

\newtheorem{rem}{Remark}[section]
\newcommand{\brem}{\begin{rem}}
\newcommand{\erem}{\end{rem}}
\newtheorem{exam}{Example}[section]

\newcommand{\bproof}{\begin{proof}}
\newcommand{\eproof}{\end{proof}}

 \title{\textbf{Generalization of soft $\mu $-compact soft generalized topological spaces}}

 \author{\textbf{Mariam Abuage $^{1}$, A. Kili\c{c}man $^{2}$, Mohammad S. Sarsak $ ^{3} $ }}
 \date{ }

\begin{document}
\maketitle
\bibliographystyle{plain}
\begin{center}
\textit{$^{1}$Institute for Mathematical Research, University Putra Malaysia, 43400 UPM, Serdang, Selangor, Malaysia} \\
\textit{$^{2}$Department of Mathematics, University Putra Malaysia, 43400 UPM, Serdang, Selangor, Malaysia}\\
\textit{$^{3}$ Department of Mathematics, The Hashemite University, P.O. Box 150459, Zarqa 13115 Jordan}
\end{center}
\medskip\noindent
\begin{abstract}
  Our work aims to introduce generalization of soft $ \mu $-compact soft generalized topological spaces, namely; soft nearly $  \mu $-compact spaces which are defined over initial universe with a fixed set of parameters. Basic properties and some significant corollaries will be introduced, Moreover, we investigated that a soft nearly $ \mu $-compact space produces a parametrized family of nearly $ \mu $-compact spaces. Some counter examples will be established to show that the converse may not be hold.
\end{abstract}

\section{Introduction}

 \ \ \ The connotation of soft set theory was first initiated by Molodtsov \cite{molodtsov1999soft} in 1999, since he proposed it as a general mathematical tool to deal with wonders while modeling the problems with incomplete informations. Therefore, a lot of researchers have been done to improve the connotation of soft set. Maji et al. \cite{maji2003soft} introduced operations of soft set. Naim Cagman et al. \cite{ccaugman2010soft} converted the definition of soft sets which is similar to that Molodtsov introduced.

 Further, general topology was expansion by several mathematical like Cs\'{a}sz\'{a}r \cite{csszar2002generalized} in 2002, who introduced the concepts of generalized topological spaces. Afterwords, a lot of authors have been done to generalize the topological notions to generalized topological setting. Sunial and Jyothis \cite{jyothis2014soft} introduced soft generalized topology on a soft set, also they  studied basic properties of soft generalized topological spaces. The same authors introduced soft $ \mu $-compact soft generalized topological spaces \cite{sunil2014soft}. In \cite{mariam2016soft} Mariam, Adem and Mohammed introduced soft $ \mu $-paracompact soft generalized topological spaces.
 Through this paper firstly, we present some substantial concepts of generalized topological spaces, which are related to our work in subsequence sections. Moreover, some primary definitions and essential results of soft set theory and soft generalized topology on an initial soft set will be given. Finally, we introduce generalization of soft $ \mu $-compact soft generalized topological spaces namely; soft nearly $ \mu $-compact spaces. Also basic properties of soft nearly $ \mu $-compact space will be introduced and some significant corollaries obtained on, Moreover, we investigated that a soft $ \mu $-paracompact space produces a parametrized family of $ \mu $-paracompact spaces.Further, subspaces of these spaces are given and some counter examples will be established to show that the converse in general may not be true.
 \section{Preliminaries}

  \ \ \ According to Cs\`{a}sz\`{a}r \cite{csszar2002generalized}, a non-empty set $ X $, $ P(X) $ denotes the power set of $ X $ and $ \mu $ be a non-empty family of $ P(X) $. The symbol $ \mu $ implies a generalized topology (briefly. GT) on $ X $, if the empty set $ \emptyset\in \mu $ and $ V_{\gamma}\in \mu $ where $ \gamma\in \Omega $ implies $ \bigcup_{\gamma\in \Omega}V_{\gamma}\in \mu $. The pair $ (X,\mu)$ is called generalized topological space (briefly. GTS) and we denote it by GTS $ (X,\mu) $ or $ X $. Each element of GT $ \mu $ is said to be $ \mu$-open set and the complement of $ \mu$-open set is called $ \mu$-closed set. Let $ A$ be a subset of a space $ (X,\mu) $, then $ i_{\mu}(A) $ (resp. $ c_{\mu}(A) $) denotes the union of all $ \mu $-open sets contained in $ A $ (resp. denotes the intersection of all $ \mu- $closed sets containing in $ A $), and $ X\backslash A $ denotes the complement of $ A $. Moreover, if a set $ X\in\mu $, then a space $ (X,\mu) $ is called $ \mu $-space \cite{noiri2006unified}. A space $(X,\mu) $ is said to be a quasi-topological space \cite{csszar2006further}, if the finite intersection of $ \mu $-open sets of $ \mu $ belongs to $ \mu $ and denoted by $ qt $-space.
   \begin{defi} \cite{sarsak2013mucompact} If a GTS $ (X,\mu)$ be a $ \mu $-space, then
   \begin{itemize}
   \item[(a)] A collection $\mathcal{V}= \{V_{\gamma}:\gamma\in\Omega\} $ of $ \mu $-open sets is called a $ \mu $-open cover of $ X$ if $ X=\bigcup_{\gamma\in\Omega}(V_{\gamma}) $.
   \item[(b)] A sub collection $\mathcal{H}= \{H_{\gamma_{i}}:\gamma_{i}\in\Omega\ for\  i\in I\} $ of a collection $ \{V_{\gamma}:\gamma\in\Omega\} $ is called $ \mu $-subcover of $ X $ if  $ X=\bigcup_{i\in I}(H_{\gamma_{i}}). $

   \end{itemize}
   \end{defi}
   \begin{defi}
    \cite{arar2014note} If a GTS $ (X,\mu) $ be a $ \mu $-space, then
       \begin{itemize}
       \item[(a)] A collection $\xi$ of subsets of $ X $ is called $ \mu $-locally finite iff for each $ x\in X $ there is $ \mu $-open set $ V $ containing $ x $  such that $ V $ intersects at finitely many elements of $ \xi. $
       \item[(b)] Let $ \mathcal{V}=\{V_{\gamma}:\gamma\in\Omega\} $  be a $ \mu $-open cover of $ X $, a collection $ \eta=\{U_{\alpha}:\alpha\in\Gamma\} $  of $ \mu $-open subsets of $ X $ is said to be a $ \mu $-open refinement of $ \mathcal{V} $ if $ \eta $ is $ \mu $-open cover of $ X $ and for each $ U\in \eta $ there is $ V\in \mathcal{V} $ such that $ U\subseteq V. $
   \end{itemize}
   \end{defi}
  We recollect the next concepts for their significance in the substance of our work, for more details see \cite{jyothis2014soft, ccaugman2010soft, maji2003soft, molodtsov1999soft}. In all the paper, $ X $ denotes initial universe,  the set of all possible parameters denoted by $ R $, $ P(X)$ is the power set of $ X $ and $ A $  be a subset of $ R $.
   \begin{defi} A soft set $ S_{A} $ on the universe $ X $ is defined by the set of ordered pairs $ S_{A}=\{(r,f_{A}(r)):r\in R, f_{A}(r)\in P(X)\} $, where $ f_{A}:R\rightarrow P(X) $ such that $ f_{A}(r)=\emptyset $ if $ r\notin A $. $ f_{A} $ is said to be approximate function of the soft set $ S_{A} $. The value of $ f_{A}(r) $ my be arbitrary. Some of them may be empty, some may be non-empty intersection. The set of all soft sets over $ X $ with $ R $ as the parameter set will be denoted by $ \mathcal{S}(X)_{R} $ or simply $ \mathcal{S}(X) $.

   \end{defi}
   \begin{defi} Let $ S_{A}\in \mathcal{S}(X) $.
   \begin{itemize}
   \item[(a)] If $ f_{A}(r)=\emptyset $ for all $ r\in R $, then $ S_{A} $ is called an empty soft set, denoted by $ S_{\emptyset} $. $ f_{A}(r)=\emptyset $ that means there is no elements in $ X $ related to the parameter $ r $ in $ R $. Moreover, we do not display such elements in the soft sets as it is meaningless to consider such parameters.
   \item[(b)] If $ f_{A}(r)=X $ for each $ r\in R $, then $ S_{A} $ is said to be an A-universal soft set, denoted by $ S_{\widehat{A}} $. If $ A=R $, then $ S_{\widehat{A}} $ is said to be an universal soft set, denoted by $ S_{\widehat{R}} $.
   \end{itemize}
   \end{defi}
   \begin{defi} If $ A\in \mathcal{S}(X) $ and $ x\in X $, then a point $ x\in S_{A} $ is said to be soft point if $ x\in f_{A}(r) $ for each $ r\in R $, and it is denoted by $ \widehat{x} $.\\

   Note: for any $ x\in X $, if $ x\notin f_{A}(r) $ for some $ r\in R $ then $ x\notin S_{A}. $

    \end{defi}
   \begin{defi} Let $ S_{A}, S_{B}\in \mathcal{S}(X) $. Then,
   \begin{itemize}
   \item[(a)] $ S_{B} $ is a soft subset of $ S_{A} $, denoted by $ S_{B}\subseteq S_{A} $, if $ f_{A}(r)\subseteq f_{B}(r) $ for all $ r\in R $.
   \item[(b)] $ S_{A}$ and $S_{B} $ are soft equal, denoted by $ S_{A}=S_{B} $, if $ f_{A}(r)=f_{B}(r) $ for all $ r\in R. $
   \item[(c)] the soft union of $ S_{A} $ and $ S_{B} $, denoted by $ S_{A}\cup S_{B} $, is defined by the approximate function $ f_{A\cup B}(r) =f_{A}(r)\cup f_{B}(r)$.
   \item[(d)] the soft intersection of $ S_{A} $ and $ S_{B} $, denoted by $ S_{A}\cap S_{B} $, is defined by the approximate function $ f_{A\cap B}(r) =f_{A}(r)\cap f_{B}(r)$.
   \item[(e)] the soft difference of $ S_{A} $ and $ S_{B} $, denoted by $ S_{A}\backslash S_{B} $, is defined by the approximate function $ f_{A\backslash B}(r) =f_{A}(r)\backslash f_{B}(r)$.
   \end{itemize}
  \end{defi}
  \begin{defi} Let $ S_{A}\in \mathcal{S}(X) $. Then,
  \begin{itemize}
  \item[(a)] the soft complement of $ S_{A} $, denoted by $ (S_{A})^{c} $, is defined by the approximate function $ f_{A^{c}}(r) =(f_{A}(r))^{c}$, where $ (f_{A}(r))^{c} $ is the complement of the set $ f_{A}(r) $, that is, $ (f_{A}(r))^{c}=X\backslash f_{A}(r) $ for all $ r\in R $.
   Clearly, $ ((S_{A})^{c})^{c}=S_{A} $ and $ (S_{\emptyset})^{c}= S_{\widehat{R}} .$
  \item[(b)] the soft power set of $ S_{A} $, denoted by $ P(S_{A}) $, is defined by $ P(S_{A})=\{S_{A_{i}}:S_{A_{i}}\subseteq S_{A}, i\in J\subseteq N\} $.

  \end{itemize}

  \end{defi}

  \begin{defi} Let $ S_{A}\in \mathcal{S}(X) $. A soft generalized topology (briefly. sGT) on $ S_{A} $, denoted by $ \mu $ or $ \mu_{S_{A}} $ is a collection of soft subsets of $ S_{A} $ such that $ S_{\emptyset}\in \mu $ and if a collection $ \{S_{A_{i}}:S_{A_{i}}\subseteq S_{A}, i\in J\subseteq N\}\subseteq\mu $ then $ \bigcup_{i\in J}(S_{A_{i}})\in \mu $.
  \end{defi}

  The pair $ (S_{A},\mu) $ is called soft generalized topological space (briefly. sGTS) and is denoted by sGTS $ (S_{A},\mu) $ or simply $ S_{A}$. It is obviously that $ S_{A}\in \mu $ must not hold.

  \begin{defi} Let a sGTS $ (S_{A},\mu) $. Then,
  \begin{itemize}
  \item[(a)] A sGT $ \mu $ is saide to be strong if $ S_{A}\in \mu$, we will denoted by strong $ \mu $-space the strong sGTS $(S_{A},\mu)$.
  \item[(b)]
  every element of $ \mu $ is called a soft $ \mu $-open set. Clearly, $ S_{\emptyset} $ is a soft $ \mu $-open set. If $ S_{B} $ be a soft subset of $ S_{A} $, then $ S_{B} $ is called soft $ \mu $-closed if its soft complement $ (S_{B})^{c} $ is a soft $ \mu $-open.

  \end{itemize}

  \end{defi}
  \begin{defi} If a sGTS $ (S_{A},\mu) $, a sub family $ \Im $ of $ \mu $ is said to be a soft basis for $ \mu $ if every member of $ \mu $ can be expressed as the soft union of some members of $ \Im $.
  \end{defi}
   \begin{defi} Let a sGTS $ (S_{A},\mu) $ and $ S_{B}\subseteq S_{A} $. Then
   \begin{itemize}
   \item[(a)] the soft union of all soft $ \mu $-open subsets of $ S_{B}$ is said to be soft $ \mu $-interior of $ S_{B} $ and denoted by $ (S_{B})^{o} $.
   \item[(b)] the soft intersection of all soft $ \mu $-closed subsets of $ S_{B} $ is said to be soft $ \mu $-closure of $ S_{B} $ and denoted by $ c(S_{B}) $.
   \item[(c)] the collection $ \mu_{S_{B}}=\{S_{V}\cap S_{B}: S_{V}\in \mu\} $ is said to be a sub-space soft generalized topology (briefly. ssGT) on $ S_{A} $. The pair $ (S_{B},\mu_{S_{B}}) $ is called soft generalized topological sub-space (briefly. sGTSS) of $ S_{A} .$
   \end{itemize}

    \end{defi}
  \begin{thm} \label{Th3.8}Let $ (S_{A},\mu) $ be a sGTS. Then a collection $\mu_{r}= \{f_{A_{i}}(r): \exists\  S_{A_{i}}\in \mu \ such\ that\ (r,f_{A_{i}}(r))\in S_{A_{i}} \ for \ i\in I\ and\ r\in R \}$
  \end{thm}

   \begin{defi} \cite{sunil2014soft} If a sGTS $ (S_{A},\mu) $ be a strong. Then
    \begin{itemize}
    \item[(a)] A collection $ \{S_{A_{i}}:i\in J\} $ of soft $ \mu $-open sets is said to be soft $ \mu $-open cover of $ S_{A} $, if $ S_{A}=\bigcup_{i\in J}(S_{A_{i}}) $.
     \item[(b)] A sub-collection $ \{S_{A_{\alpha_{i}}}:\alpha_{i}\in \Delta\ for\ i\in J\} $ of a collection $ \{S_{A_{\alpha}}:\alpha\in \Delta\} $ is saide to be soft $ \mu $-open sub-cover of $ S_{A} $, if $ S_{A}=\bigcup_{i\in J}\{S_{A_{\alpha_{i}}}: \alpha_{i}\in \Delta\}.$
     \item[(c)] A sGTS $ (S_{A},\mu) $ is said to be soft $ \mu $-compact on $ S_{A} $ if for each soft $ \mu $-open cover of $ S_{A} $ admits a finite soft $ \mu $-open sub-cover.
    \end{itemize}
    \end{defi}

    \begin{thm} \label{Th2.14} \cite{mariam2016soft} Let a strong sGTS $ (S_{\widehat{A}},\mu) $, then \begin{itemize}
   \item[(a)]
     a collection $ S_{A_{i}
       }=\{(r,f_{A_{i}}(r)):i\in J\} $ for each $ r\in R $ is a soft subsets of $ S_{\widehat{A}} $ if and only if a collection $ \{f_{A_{i}}(r):i\in J\}$ for each $ r\in R  $, is subsets of a $ \mu_{r} $-space $(X,\mu_{r}). $

        \item[ (b)] a collection $ S_{A_{i}
    }=\{(r,f_{A_{i}}(r)):i\in J\} $ for each $ r\in R $ is a soft $ \mu $-open cover of $ S_{\widehat{A}} $ if and only if a collection $ \{f_{A_{i}}(r):i\in J\}$ for each $ r\in R  $, is $ \mu_{r} $-open cover of a $ \mu_{r} $-space $(X,\mu_{r}). $
    \end{itemize}
     \end{thm}

\section{Soft $ n\mu $-compact spaces and subspaces}

  \ \ \ To simply , sGTS $ (S_{A},\mu) $ always means strong sGTS. if not we cite for that.
 \begin{defi} A soft subset $ S_{B} $ of a sGTS $ (S_{A},\mu) $ is said to be soft $ \mu $-regular open (resp. soft $ \mu $-regular closed) iff $ S_{B}= (c(S_{B}))^{o} $ (resp. $ S_{B}= c((S_{B})^{o} )$)

 \end{defi}

Obviously, every soft $ \mu $-regular open (resp. soft $ \mu $-regular closed) is soft $ \mu $-open (resp. soft $ \mu $-closed) set, the following example shows that the converse not be hold:
\begin{exam} Let $ X=\{a,b,c\}, R=\{r_{1},r_{2},r_{3}\}, A=\{r_{1},r_{2}\}\subseteq R $ and $ S_{A} =\{(r_{1},X),(r_{2},\{b,c\})\}$. Let $ \mu=\{S_{\emptyset},S_{A_{1}},S_{A_{2}},S_{A_{3}},S_{A}\} $ where \\ $ S_{A_{1}}=\{(r_{1},\{b\}),(r_{2},\{b,c\})\}, $\\
$ S_{A_{2}}=\{(r_{1},\{a,c\}),(r_{2},\{c\})\},$\\
$ S_{A_{3}}=\{(r_{1},\{a,b\}),(r_{2},\{b,c\})\}.$\\

Then the pair $ (S_{A},\mu) $ is sGTS over $ X $, and it clear that $ S_{A_{1}},S_{A_{2}}, S_{A_{3}}  $ are soft $ \mu $-open sets but not soft $ \mu $-regular open.
\end{exam}
\begin{thm} Let a sGTS $ (S_{A},\mu) $, if \begin{itemize}
\item[(i)] $ S_{F} $ is soft $ \mu $-closed set then $ (S_{F})^{o} $ is soft $ \mu $-regular open.
\item[(ii)] $ S_{V} $ is soft $ \mu $-open set then $ c(S_{V}) $ is soft $ \mu $-regular closed.

\end{itemize}
\end{thm}
{\bf Proof.} $ (i) $ Define $ S_{B}= (S_{F})^{o}$. $ S_{B}\subseteq S_{F}$ implies $ S_{B}\subseteq c(S_{B})\subseteq S_{F}$ and, by $ S_{B}\in \mu $, $ S_{B}\subseteq (c(S_{B}))^{o}\subseteq (S_{F})^{o}=S_{B} $. Thus $ S_{B}= (c(S_{B}))^{o}$. \\

$ (ii) $ Define $ S_{F}= c(S_{V})$. $ S_{V}\subseteq S_{F}$ implies $ S_{V}\subseteq (S_{F})^{o}\subseteq S_{F}$ and $ S_{F}= c(S_{V})\subseteq c((S_{F})^{o})\subseteq c(S_{F})=S_{F} $ as $ S_{F} $ is soft $ \mu $-closed. Thus $ S_{F}= c((S_{F})^{o})$. \\

If $ A\in \mathcal{S}(X) $ be A-universal soft set, then, on using Theorem \ref{Th2.14} (a) we conclude the corollary below:
\begin{corr} \label{cor3.1}
Let a sGTS $ (S_{\widehat{A}},\mu) $, then a collection $ S_{A_{i}
       }=\{(r,f_{A_{i}}(r)):i\in J\} $ for each $ r\in R $ is a soft $ \mu $-regular sets of $ S_{\widehat{A}} $ if and only if a collection $ \{f_{A_{i}}(r):i\in J\}$ for each $ r\in R  $, is $ \mu_{r} $-regular sets of a $ \mu_{r} $-space $(X,\mu_{r}) $ for each $ r\in R. $
\end{corr}

 \begin{defi} A sGTS $ (S_{A},\mu) $ is called soft nearly $ \mu $-compact (briefly. soft $ n\mu $-compact) if for each soft $ \mu $-open cover $ \{S_{A_{i}}:i\in J\subseteq \mathcal{N}\} $ of $ S_{A} $ admits a finite soft $ \mu $-open sub-cover $ \{S_{A_{k}}:k=1,2,.....,n\} $ such that $$ S_{A}=\bigcup_{k=1}^{n}(c(S_{A_{k}}))^{o}. $$
 i.e. for each soft $ \mu $-regular open cover of $ S_{A} $ admits a finite soft $ \mu $-open sub-cover.
 \end{defi}

 It is easily to note that every soft $ \mu $-compact space is soft $ n\mu $-compact, the converse in general is not true as in the example shows:
  \begin{exam} Let $ X=\mathcal{N}, R= A =\{r_{i}:i\in \mathcal{N}\}$ and $ S_{\widehat{R}}=\{(r_{i},X):r_{i}\in R\} $, let $ \Im=\{(r,\{1,x\}):x\in X, x\neq 1\} $ for each $ r\in R $. Consider a sGT $ \mu(\Im) $ generated on sGTS $ S_{\widehat{R}} $ by the soft basis $ \Im $. Then, only $ S_{\widehat{R}} $ and $ S_{\emptyset} $ are soft $ \mu $-regular open set so a sGTS $ (S_{\widehat{R}},\mu(\Im)) $ is soft $ n\mu(\Im) $-compact but it is not soft $ \mu(\Im) $-compact, since a collection $ \{S_{\widehat{R}_{i}}:i\in \mathcal{N}\} $, where \\$ S_{\widehat{R}_{1}}=\{(r_{1},\{1,2\}),(r_{2},\{1,2,3\}),(r_{3},\{1,2,3,4\}),......\} $,\\
  $ S_{\widehat{R}_{2}}=\{(r_{1},\{1,3\}),(r_{2},\{1,2,4\}),(r_{3},\{1,2,3,5\}),......\} ,$\\
  $S_{\widehat{R}_{3}}=\{(r_{1},\{1,4\}),(r_{2},\{1,2,5\}),(r_{3},\{1,2,3,6\}),......\} ,$\\
  .\\
  .\\
  .\\
  is soft $ \mu(\Im)$-open cover of sGTS $ (S_{\widehat{R}},\mu(\Im))$ with no finite soft $ \mu(\Im) $-open sub-cover.
  \end{exam}

   \begin{thm} A  sGTS $ (S_{A},\mu) $ is soft $ n\mu $-compact $\Leftrightarrow $ for every collection of soft $ \mu$-regular closed subsets of $ S_{A} $ having the soft finite intersection property, has a nonempty soft intersection.
   \end{thm}
   {\bf proof.}   $ (\Rightarrow) $ Let $ \{S_{F_{h}}:h\in J\subseteq \mathcal{N}\} $ is a collection of soft $\mu $-regular closed sets of $ S_{A} $ and suppose that $ \bigcap_{h\in J}(S_{F_{h}})= S_{\emptyset} $. Thus $ \{S_{A}\backslash S_{F_{h}}:h\in J\} $ is a soft $ \mu $-regular open cover of the soft $ n\mu $-compact space $ S_{A} $, so it has a finite soft $ \mu $-sub-cover $ \{S_{A}\backslash S_{F_{k}}:k=1,2,...,n\} $ such that $ \bigcup_{k=1}^{n}(S_{A}\backslash S_{F_{k}})=S_{A} $. This implies that $ \bigcap_{k=1}^{n}(S_{F_{k}})=S_{A}\backslash S_{A}=S_{\emptyset} $, and this is a contradiction with the soft finite intersection property. Then $ \bigcap_{h\in J}(S_{F_{h}})\neq S_{\emptyset} $.\\

   $ (\Leftarrow) $ Let that a sGTS $ (S_{A},\mu) $ is not soft $ n\mu $-compact, then there exists a soft $ \mu $-regular open cover $ \{S_{V_{i}}:i\in J\subseteq \mathcal{N}\} $ with no finite soft $ \mu$-sub-cover. So $ S_{A}\nsubseteq \bigcup_{k=1}^{n}(S_{V_{k}}) $ for any countable soft subset $ \{S_{V_{k}}:k\in \mathbb{N}\} $ of $ \{S_{V_{i}}:i\in J\} $. then it follows that $$ S_{A}\backslash \bigcup_{k=1}^{n}(S_{V_{k}})=\bigcap_{k=1}^{n}(S_{A}\backslash S_{V_{k}})\neq S_{\emptyset} .$$ Thus $ \{S_{A}\backslash S_{V_{i}}:i\in J\} $ is a collection of soft $ \mu $-regular closed sets of $ S_{A}$ satisfies the soft finite intersection property. So, by hypothesis $ \bigcap_{i\in J}(S_{A}\backslash S_{V_{i}})\neq S_{\emptyset} $. It means that $ S_{A}\backslash(\bigcup_{i\in J}(S_{V_{i}}))\neq S_{\emptyset} $, which is a contradiction with the fact that $ \{S_{V_{i}}:i\in J\} $ is a soft $ \mu$-cover of $ S_{A}$. This implies that a sGTS $ S_{A} $ is soft $ n\mu $-compact.
   \begin{defi} A $ \mu $-space $ (X,\mu) $ is said to be
   \begin{itemize}
   \item[(a)] nearly $ \mu $-compact (briefly. $ n\mu $-compact) if for each $ \mu $-open cover $ \{V_{\gamma}:\gamma\in\Omega\} $ of $ X $ admits a finite $ \mu $-subcover $ \{V_{\gamma_{k}}: k=1,2,....n\} $ such that $$ X=\bigcup_{k=1}^{n}(i_{\mu}(c_{\mu}(V_{\gamma_{k}}))).$$
   \item[(b)] nearly $ \mu $-paracompact space (briefly. $ n\mu $-paracompact) if for every $ \mu $-regular open cover of $ X $ admits a $ \mu $-open $ \mu $-locally refinement.

   \end{itemize}
   \end{defi}
   \begin{rem}
     It is clear to see that every $ n\mu $-compact space is $ n\mu $-paracompact, because any finite $ \mu$-sub-cover is $ \mu $-open $ \mu $-locally finite refinement, in general the inversion not be true as an example below:
      \end{rem}
      \begin{exam} If $ X= \mathbb{N} $ and $\beta =\{\{n,n+1\}:n\in \mathbb{N}\}$. Consider a GT $ \mu(\beta) $ generated on $ X $ by the $ \mu $-base $ \beta $, then a space $ (X,\mu(\beta)) $ is $ n\mu $-paracompact but it is not $ n\mu $-compact, since $ \{\{n,n+1\}: n\ is\ an\ odd\  number\} $ is  $ \mu $-regular open cover of $ X $ with no finite $ \mu $-sub-cover.
      \end{exam}

      Manner of proof, of the next Theorem is similar to proof Theorem...[],so omitted.
      \begin{thm} Let a sGTS $ (S_{A},\mu) $ be a soft $ n\mu $-compact, then a GTS $ (X,\mu_{r}) $ is $ n\mu_{r} $-compact (resp. $ n\mu_{r} $-paracompact) for each $ r\in R. $
       \end{thm}

       In Theorem above, we observe that corresponding to each parameter $ r\in R $, we have a $ n\mu_{r} $-compact (resp. $ n\mu_{r} $-paracompact) GTS. Then a sGTS $ (S_{A},\mu) $ gives a parametrized family of $ n\mu_{r} $-compact (resp. $ n\mu_{r} $-paracompact) GTS.

      On using Theorem \ref{Th2.14} (b), Corollary \ref{cor3.1} and the remark above, we obtained on the corollary below:
      \begin{corr} Let  $ A\in \mathcal{S}(X)$ be A-universal soft set, then a sGTS $ (S_{\widehat{A}},\mu) $ is a soft $ n\mu $-compact $ \Leftrightarrow $ a $ \mu_{r} $-space $ (X,\mu_{r}) $ is a $ n\mu_{r} $-compact (resp. $ n\mu_{r} $-paracompact) for each $ r\in R .$

      \end{corr}
    \begin{defi}  Let $ S_{B} $ be a soft subset of sGTS $ (S_{A},\mu) $. A sGTSS $
        (S_{B},\mu_{S_{B}}) $ of a sGTS $ (S_{A},\mu) $ is called soft $ \mu_{S_{B}} $-compact (resp. soft $ n\mu_{S_{B}} $-compact) if for any soft $ \mu_{S_{B}} $-open (resp. soft $ \mu_{S_{B}} $-regular open) cover of $ S_{B} $ admits a finite soft $ \mu_{S_{B}} $-open sub-cover.
      \end{defi}

   \begin{rem} In \cite{sunil2014soft}, it was pointed out that a soft $ \mu $-closed subset of a soft $ \mu $-compact space is soft $ \mu $-compact but in soft $ n\mu $-compact spaces it is not necessary true, as the example below:
   \end{rem}
   \begin{exam} Let $ X=\mathbb{R}$ be the set of real numbers, $ A=R=\{r_{i}:i\in \mathcal{N}\} $ and $ S_{\widehat{R}}=\{(r_{i},X): r_{i}\in R\} $, let $ \Im =\{(r,A):A\subseteq X,A^{c}\ is\ a\ countable\ soft\ set\}$ for each $ r\in R $. Consider a sGT $ \mu(\Im) $ generated on sGTS $ (S_{\widehat{R}},\mu(\Im)) $ by the soft basis $ \Im $. Clearly, a sGTS $ (S_{\widehat{R}},\mu(\Im)) $ is soft $ n\mu $-compact space, since the only soft $ \mu $-regular open cover is the set $ S_{\widehat{R}} $ itself, which has only one element so it is soft finite. Now, let $ S_{F}= \{(r_{i},\mathcal{N}): r_{i}\in R\} $, then $ S_{F} $ is soft $ \mu $-closed but not soft $ \mu $-open set and let the induced soft generalized topology on $ S_{F} $ is generated by by the soft basis $ \hbar=\{(r,\{x\}):x\in \mathcal{N}\} $ for each $ r\in R $. Then a sGTSS  $ (S_{F},\mu_{S_{F}}(\hbar)) $ is not soft $ n\mu_{S_{F}}(\hbar) $-compact, since a collection $ \{S_{F_{i}}:i\in \mathcal{N}\} $, where\\
   $ S_{F_{1}}=\{(r_{1},\{1\}),(r_{2},\{1,2\}),(r_{3},\{1,2,3\}),....\} $\\
   $S_{F_{2}}=\{(r_{1},\{2\}),(r_{2},\{1,3\}),(r_{3},\{1,2,4\}),....\} $ \\
   $S_{F_{3}}=\{(r_{1},\{3\}),(r_{2},\{1,4\}),(r_{3},\{1,2,5\}),....\} $\\
   .\\
   .\\
   .\\
   is a soft $ \mu_{S_{F}}(\hbar) $-open cover of $ S_{F} $ with no finite soft $ \mu_{S_{F}}(\hbar) $-open sub-cover.
   \end{exam}
   \begin{defi} A soft subset $ S_{B} $ of a sGTS $ (S_{A},\mu) $ is called soft $ \mu $-clopen if it is both soft $ \mu $-open and soft $ \mu $-closed.

   \end{defi}
   \begin{thm} Let a sGTS $ (S_{A},\mu) $ and a set $ S_{B} $ be a soft $ \mu $-open subset of $ S_{A}$. Then soft $ \mu_{S_{B}} $-regular open (resp. soft $ \mu_{S_{B}} $-regular closed) sets in the induced soft generalized topology $ (S_{B},\mu_{S_{B}}) $ are the form $ S_{B}\cap S_{U} $, where $ S_{U} $ is soft $ \mu $-regular open (resp. soft $ \mu $-regular closed) in sGTS $ (S_{A},\mu) $.

   \end{thm}
   {\bf Proof.} Let $ S_{V} $ be a soft $ \mu_{S_{B}} $-regular open set then $ c'(S_{V})=c(S_{V})\cap S_{B} $, hence $ (c'(S_{V}))^{o'}=(c(S_{V}))^{o}\cap S_{B}= S_{U}\cap S_{B} $, since $ c'(S_{V})$ (resp. $(S_{V})^{o'} $) denotes the soft $ \mu_{S_{B}} $-clouser (resp. soft $ \mu_{S_{B}} $-interior) set in sGTSS $ (S_{B},\mu_{S_{B}}) $ and $ S_{U}=(c(S_{V}))^{o} $ is soft $ \mu $-regular open set. Same proof of soft $ \mu_{S_{B}} $-regular closed set.\\
   Conversely, let $ S_{U} $ be a soft $ \mu $-regular open set and $ S_{V}=S_{U}\cap S_{B} $. Then $$ (c'(S_{V}))^{o'}=(c(S_{U}\cap S_{B})\cap S_{B})^{o}=(c(S_{U})\cap S_{B})^{o}=(c(S_{U}))^{o}\cap S_{B}=S_{U}\cap S_{B}=S_{V} ,$$ hence $ S_{V} $ is soft $ \mu_{S_{B}} $-regular open set. Same proof of soft $ \mu $-regular closed set.\\

   \ \ \ The following theorem gives an analogous result for soft $ n\mu $-compact subspaces.
   \begin{thm} Let a sGTS $ (S_{A},\mu) $ be a soft $ n\mu $-compact and a set $ S_{B} $ be a soft $ \mu $-clopen subset of $ S_{A}$. Then a soft generalized topological space $ (S_{B},\mu_{S_{B}})$ is  soft $n \mu_{S_{B}} $-compact subspace.
   \end{thm}
   {\bf Proof.} Let $ \{S_{V_{i}}:i\in J\subseteq \mathcal{N}\} $ be a soft $ \mu_{S_{B}}$-regular open cover of $S_{B}$. By Theorem ..., for each $ i\in J$,  $S_{V_{i}}=S_{U_{h}}\cap S_{B}$ where $ S_{U_{h}} $ is soft $ \mu $-regular open sets in $ S_{A} $ for each $ h\in J $ . Thus the family $ \{S_{A}\backslash S_{B}\}\cup\{S_{U_{h}}:h\in J\} $ forms a soft $ \mu $-regular open cover of the soft $ n\mu $-compact space $ S_{A} $. So there is a finite soft $ \mu $-sub-cover $ \{S_{A}\backslash S_{B}\}\cup\{S_{U_{k}}:k=1,2,...,n\} $. Thus $ S_{B}=\bigcup_{k=1}^{n}(S_{U_{k}}\cap S_{B})=\bigcup_{k=1}^{n}(S_{V_{k}}) $. It follows that a sGTSS $ (S_{B},\mu_{S_{B}}) $ is soft $ n\mu_{S_{B}} $-compact.
   \section{Conclusion}
     \ \ \ In this study, firstly, we present some substantial concepts of generalized topological spaces, which are related to our work. Moreover, some primary definitions and essential results of soft set theory and soft generalized topology on an initial soft set are given. Secondly, we introduced soft nearly $ \mu $-compact soft generalized topological spaces and investigated that a soft $ \mu $-paracompact space produces a parametrized family of $ \mu $-paracompact spaces. Finally, subspaces of these spaces are given and some counter examples established to show that the converse may not be hold.
%

\end{document}